\documentclass[english]{article}
\usepackage[T1]{fontenc}
\usepackage[latin1]{inputenc}
\usepackage{float}
\usepackage{graphicx}
\usepackage{setspace}
\usepackage{amssymb}

\makeatletter

\newcommand{\lyxline}[1]{
  {#1 \vspace{1ex} \hrule width \columnwidth \vspace{1ex}}
}

\providecommand{\tabularnewline}{\\}

\newcommand{\lyxaddress}[1]{
\par {\raggedright #1
\vspace{1.4em}
\noindent\par}
}

\usepackage{babel}
\makeatother
\begin{document}

\title{About Goldbach strong conjecture}

\author{\textit{G. Funes}$^{A,B}$\textit{, D. Gulich}$^{A,B}$\textit{,
L. Garavaglia}$^{C}$\textit{, y M. Garavaglia}$^{A,B}$}

\maketitle
\begin{singlespace}
\begin{center}
{\scriptsize $^{A}$Departamento de Física, Facultad de Ciencias Exactas,
Universidad Nacional de La Plata, Argentina}
\par\end{center}{\scriptsize \par}

\begin{center}
{\scriptsize $^{B}$Laboratorio de procesamiento Láser, Centro de
investigaciones Ópticas, La Plata, Argentina}
\par\end{center}{\scriptsize \par}

\begin{center}
{\scriptsize $^{C}$Aranjuez, España}{\footnotesize }\\

\par\end{center}{\footnotesize \par}
\end{singlespace}

\lyxaddress{\emph{\footnotesize E-mail addresses:} {\footnotesize Gustavo Funes:
gfunes@ciop.unlp.edu.ar; Damián Gulich: dgulich@ciop.unlp.edu.ar;
Mario Garavaglia: garavagliam@ciop.unlp.edu.ar}}

\lyxline{\normalsize}

\begin{abstract}
In this work we use the number classification in families of the form
$6n+1$, and $6n+5$ with $n$ integer (Such families contain all
odd prime numbers greater than 3 and other compound numbers related
with primes). We will use this kind of classification in order to
attempt an approach to Goldbach strong conjecture. By means of a geometric
method of binary bands of numbers we conceive a new form of study
of the stated problem. 
\end{abstract}
\lyxline{\normalsize}

\section{Introduction}

In June 1724, the mathematician Christian Goldbach wrote a letter
to Leonhard Euler in which he proposed the following conjecture:

\begin{verse}
{}``All odd integer greater than two can be expressed as the sum
of three primes''
\end{verse}
The relation proposed by Goldbach is based on the primality of number
$1$, which was considered a prime number at those times. Nevertheless,
this was rejected afterwards. 

A modern version of Goldbach original conjecture is:

\begin{verse}
{}``All integers greater than 5 can be expressed as the sum of three
primes'',
\end{verse}
which is also called {}``Triple'' or {}``Ternary'' conjecture.

Thus, Euler proposed that this version has an origin on a strong conjecture:

\begin{verse}
{}``All even numbers greater than two can be written as the sum of
two primes'',
\end{verse}
which is called {}``Double'' or {}``Binary'' conjecture.

It is clear that the {}``ternary'' conjecture is a consequence of
the {}``double'' one, for that reason the last one is called {}``strong''
and the ternary is called {}``weak''.

\section{About the conjecture}

We will begin our analysis of the strong Goldbach conjecture by using
the results stated by Garavaglia \textit{et al} \cite{Garavaglia},
whom had established that all numbers can be classified in families
or kinds, which obey the following formulas:

\begin{equation}
\alpha=1+6n\label{eq:alfa}\end{equation}

\begin{equation}
\epsilon=2+6n\label{eq:epsilon}\end{equation}

\begin{equation}
\gamma=3+6n\label{eq:gamma}\end{equation}

\begin{equation}
\delta=4+6n\label{eq:delta}\end{equation}

\begin{equation}
\beta=5+6n\label{eq:beta}\end{equation}

\begin{equation}
\zeta=6+6n\label{eq:zeta}\end{equation}

All odd integers are represented by equations \ref{eq:alfa}, \ref{eq:gamma},
and \ref{eq:beta}, while even integers do by \ref{eq:epsilon}, \ref{eq:delta},
and \ref{eq:zeta}.

Since only families $\alpha$ and $\beta$ contain all prime numbers
except $+2$ which is $\epsilon(n=0)$, $-2$ who belongs to $\delta(n=-1)$,
$+3$ which is $\gamma(n=0)$ and $-3$ who is also $\gamma(n=-1)$
\cite{Garavaglia}. Using the definitions given by the formulas \ref{eq:alfa}
and \ref{eq:beta}, we will use in the following:

\[
2k=p_{1}+p_{2}\]

\noindent where $k$ is an integer greater then one and $p_{1,2}$
are two prime numbers. Taking into account that most primes belong
to the families $\alpha$ and $\beta$, the possibilities of sum will
be the following:

\[
2k=\alpha+\alpha\]

\[
2k=\alpha+\beta\]

\[
2k=\beta+\beta\]

We will replace the expressions of $2k$ in the formulas \ref{eq:alfa}
and \ref{eq:beta}, to obtain:

\[
2k=\alpha+\alpha=(6n_{1}+1)+(6n_{2}+1)\]

\[
2k=\alpha+\beta=(6n_{1}+1)+(6n_{2}+5)\]

\[
2k=\beta+\beta=(6n_{2}+5)+(6n_{2}+5)\]

\noindent rearranging:

\[
6(n_{1}+n_{2})+2=2\times(3(n_{1}+n_{2})+1)\]

\[
6(n_{1}+n_{2})+6=2\times(3(n_{1}+n_{2})+3)\]

\[
6(n_{1}+n_{2})+10=2\times(3(n_{1}+n_{2})+5)\]

\noindent It is demonstrated that the strong conjecture is compatible
with the classification in families, it means that a sum of two numbers,
$\alpha$ o $\beta$ results of the form $2\times k$, with $k$ integer.
The preceding formulas (not counting $2\times$) should generate all
integers $k$ with $n_{1}$ and $n_{2}$ integers.

\section{About number generation}

We will define:

\begin{equation}
\vartheta_{1}=(3(n_{1}+n_{2})+1)\label{eq:generador 1}\end{equation}

\[
\vartheta_{2}=(3(n_{1}+n_{2})+3)\]

\[
\vartheta_{3}=(3(n_{1}+n_{2})+5)\]

\noindent which justify the following:

\begin{equation}
\vartheta_{2}=\vartheta_{1}+2\label{eq:generador 2}\end{equation}

\begin{equation}
\vartheta_{3}=\vartheta_{1}+4\label{eq:generador 3}\end{equation}

\noindent in order to maintain a constant increment we could replace
the last equation by:

\[
\vartheta_{3}=\vartheta_{2}+2\]

In Table \ref{tab:Numeros-generados-por} we will use equations \ref{eq:generador 1},
\ref{eq:generador 2}, and \ref{eq:generador 3} to obtain the numbers
which are generated by each one.

\noindent \begin{center}
\begin{table}[h]
\noindent \begin{centering}
\begin{tabular}{|c|c|c|c|}
\hline 
$n_{1}$/ $n_{2}$&
$\vartheta_{1}$&
$\vartheta_{2}$&
$\vartheta_{3}$\tabularnewline
\hline
\hline 
0 / 0&
1&
3&
5\tabularnewline
\hline 
1 / 0&
4&
6&
8\tabularnewline
\hline 
1 / 1&
7&
9&
11\tabularnewline
\hline 
1 / 2&
10&
12&
14\tabularnewline
\hline 
2 / 2&
13&
15&
17\tabularnewline
\hline 
2 / 3&
16&
18&
20\tabularnewline
\hline 
3 / 3&
19&
21&
23\tabularnewline
\hline 
3 / 4&
22&
24&
26\tabularnewline
\hline 
4 / 4&
25&
27&
29\tabularnewline
\hline 
5 / 4&
28&
30&
32\tabularnewline
\hline 
5 / 5&
31&
33&
35\tabularnewline
\hline 
5 / 6&
34&
36&
38\tabularnewline
\hline 
6 / 6&
37&
39&
41\tabularnewline
\hline 
6 / 7&
40&
42&
44\tabularnewline
\hline 
7 / 7&
43&
45&
47\tabularnewline
\hline 
7 / 8&
46&
48&
50\tabularnewline
\hline 
8 / 8&
49&
51&
53\tabularnewline
\hline 
8 / 9&
52&
54&
56\tabularnewline
\hline 
9 / 9&
55&
57&
59\tabularnewline
\hline
\end{tabular}
\par\end{centering}

\caption{\label{tab:Numeros-generados-por}Generated numbers by $\vartheta$
equations.}
\end{table}

\par\end{center}

Then formulas \ref{eq:generador 1}, \ref{eq:generador 2}, and \ref{eq:generador 3},
have generated all numbers except $2$. In order to obtain it, and
to obtain all negative numbers as well, we will add negative $n_{i}$
to Table \ref{tab:Numeros-generados-por}. By doing this we obtain
the following results indicated in Table \ref{tab:Numeros-generados-con negativos}.

\begin{center}
\begin{table}[H]
\noindent \begin{centering}
\begin{tabular}{|c|c|c|c|}
\hline 
$n_{1}$/ $n_{2}$&
$\vartheta_{1}$&
$\vartheta_{2}$&
$\vartheta_{3}$\tabularnewline
\hline
\hline 
-1 / -2&
-8&
-6&
-4\tabularnewline
\hline 
-1 / -1&
-5&
-3&
-1\tabularnewline
\hline 
-1 / 0&
-2&
0&
2\tabularnewline
\hline 
0 / 0&
1&
3&
5\tabularnewline
\hline 
1 / 0&
4&
6&
8\tabularnewline
\hline
1 / 1&
7&
9&
11\tabularnewline
\hline
\end{tabular}
\par\end{centering}

\caption{\label{tab:Numeros-generados-con negativos}Generated numbers by
$\vartheta$ equations and negative $n_{i}$ }
\end{table}

\par\end{center}

Taking into account that differences between numbers expressed by
equations $\vartheta_{i}$ is always 2, according to the formulas
\ref{eq:generador 1}, \ref{eq:generador 2}, and \ref{eq:generador 3},
and that the differences between numbers of the same family is always
3, $\vartheta_{1}(1)-\vartheta_{1}(2)=3$. We conclude that these
generator families ($\vartheta_{i}$) will always cover the entire
coordinate axis of integer numbers (Figure \ref{fig:Cobertura-del-eje}).

\noindent \begin{center}
\begin{figure}[h]
\noindent \begin{centering}
\includegraphics[width=8cm,keepaspectratio]{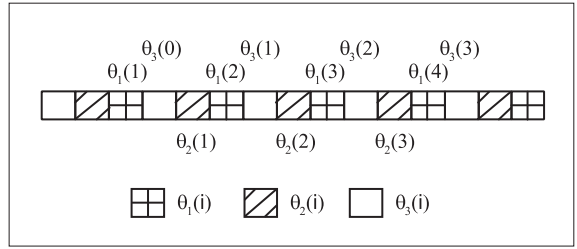}
\par\end{centering}

\caption{\label{fig:Cobertura-del-eje} Axis cover of $\mathbb{Z}$ numbers
by means of $\vartheta_{i}$ generators }
\end{figure}

\par\end{center}

Then, with the introduction of the classification of odd integer numbers
like $\alpha$ and $\beta$, we have obtained a result that seems
valid to begin analyze the Goldbach strong conjecture.

\section{About the restriction of our statements}

It is true that we solve the problem of obtaining even numbers by
adding number families which contain almost all prime numbers ($\alpha$
and $\beta$). Nevertheless, we have skipped an important point, not
all $\alpha$ and $\beta$ are primes. It is true that these families
contain all prime numbers except $\pm2$ and $\pm3$, but they also
contain such compound numbers which are product of elements of these
very families \cite{Garavaglia}, ($\alpha\times\alpha$, $\beta\times\beta$,
and $\alpha\times\beta$). Actually we have demonstrated that all
even numbers can be expressed by the sum of two $\alpha$ , two $\beta$
, or a $\beta$ plus an $\alpha$. Eventually, the problem of this
new form to analyze Goldbach strong conjecture, is based on the possibility
of replacing compound numbers by prime ones while keeping the sum
intact. This fact could not be estimated \textit{a priori} , because
not all compound numbers are followed or preceded by primes. In other
words, what we try to do is to maintain certain sums avoiding the
use of compound elements from $\alpha$ or $\beta$. We will see that
it is possible to apply this proceeding to some examples, but our
primary objective is to prove that it is valid to all sums of odd
numbers. 

\begin{table}[h]
\begin{centering}
\begin{tabular}{|c|c|c|c|c|}
\hline 
{\footnotesize $\beta+\beta$}&
{\footnotesize $\beta_{1}$}&
{\footnotesize $\beta_{2}$}&
&
{\footnotesize $\beta_{primo}+\beta_{primo}$}\tabularnewline
\hline
\hline 
{\footnotesize 10}&
{\footnotesize 5}&
{\footnotesize 5}&
&
\tabularnewline
\hline 
{\footnotesize 16}&
{\footnotesize 5}&
{\footnotesize 11}&
&
\tabularnewline
\hline 
{\footnotesize 22}&
{\footnotesize 11}&
{\footnotesize 11}&
&
\tabularnewline
\hline 
{\footnotesize 28}&
{\footnotesize 11}&
{\footnotesize 17}&
&
\tabularnewline
\hline 
{\footnotesize 34}&
{\footnotesize 17}&
{\footnotesize 17}&
&
\tabularnewline
\hline 
{\footnotesize 40}&
{\footnotesize 17}&
{\footnotesize 23}&
&
\tabularnewline
\hline 
{\footnotesize 46}&
{\footnotesize 23}&
{\footnotesize 23}&
&
\tabularnewline
\hline 
{\footnotesize 52}&
{\footnotesize 23}&
{\footnotesize 29}&
&
\tabularnewline
\hline 
{\footnotesize 58}&
{\footnotesize 29}&
{\footnotesize 29}&
&
\tabularnewline
\hline 
{\footnotesize 64}&
{\footnotesize 29}&
{\footnotesize 35}&
{\footnotesize 64}&
{\footnotesize 23+41}\tabularnewline
\hline 
{\footnotesize 70}&
{\footnotesize 35}&
{\footnotesize 35}&
{\footnotesize 70}&
{\footnotesize 29+41}\tabularnewline
\hline 
{\footnotesize 76}&
{\footnotesize 35}&
{\footnotesize 41}&
{\footnotesize 76}&
{\footnotesize 29+47}\tabularnewline
\hline 
{\footnotesize 82}&
{\footnotesize 41}&
{\footnotesize 41}&
&
\tabularnewline
\hline 
{\footnotesize 88}&
{\footnotesize 41}&
{\footnotesize 47}&
&
\tabularnewline
\hline 
{\footnotesize 94}&
{\footnotesize 47}&
{\footnotesize 47}&
&
\tabularnewline
\hline 
{\footnotesize 100}&
{\footnotesize 47}&
{\footnotesize 53}&
&
\tabularnewline
\hline 
{\footnotesize 106}&
{\footnotesize 53}&
{\footnotesize 53}&
&
\tabularnewline
\hline 
{\footnotesize 112}&
{\footnotesize 53}&
{\footnotesize 59}&
&
\tabularnewline
\hline 
{\footnotesize 118}&
{\footnotesize 59}&
{\footnotesize 59}&
&
\tabularnewline
\hline 
{\footnotesize 124}&
{\footnotesize 59}&
{\footnotesize 65}&
{\footnotesize 124}&
{\footnotesize 53+71}\tabularnewline
\hline 
{\footnotesize 130}&
{\footnotesize 65}&
{\footnotesize 65}&
{\footnotesize 130}&
{\footnotesize 59+71}\tabularnewline
\hline 
{\footnotesize 136}&
{\footnotesize 65}&
{\footnotesize 71}&
{\footnotesize 136}&
{\footnotesize 53+83}\tabularnewline
\hline 
{\footnotesize 142}&
{\footnotesize 71}&
{\footnotesize 71}&
&
\tabularnewline
\hline 
{\footnotesize 148}&
{\footnotesize 71}&
{\footnotesize 77}&
{\footnotesize 148}&
{\footnotesize 59+89}\tabularnewline
\hline 
{\footnotesize 154}&
{\footnotesize 77}&
{\footnotesize 77}&
{\footnotesize 154}&
{\footnotesize 71+83}\tabularnewline
\hline 
{\footnotesize 160}&
{\footnotesize 77}&
{\footnotesize 83}&
{\footnotesize 160}&
{\footnotesize 71+89}\tabularnewline
\hline 
{\footnotesize 166}&
{\footnotesize 83}&
{\footnotesize 83}&
&
\tabularnewline
\hline 
{\footnotesize 172}&
{\footnotesize 83}&
{\footnotesize 89}&
&
\tabularnewline
\hline 
{\footnotesize 178}&
{\footnotesize 89}&
{\footnotesize 89}&
&
\tabularnewline
\hline 
{\footnotesize 184}&
{\footnotesize 89}&
{\footnotesize 95}&
{\footnotesize 184}&
{\footnotesize 83+101}\tabularnewline
\hline 
{\footnotesize 190}&
{\footnotesize 95}&
{\footnotesize 95}&
{\footnotesize 190}&
{\footnotesize 89+101}\tabularnewline
\hline 
{\footnotesize 196}&
{\footnotesize 95}&
{\footnotesize 101}&
{\footnotesize 196}&
{\footnotesize 89+107}\tabularnewline
\hline 
{\footnotesize 202}&
{\footnotesize 101}&
{\footnotesize 101}&
&
\tabularnewline
\hline 
{\footnotesize 208}&
{\footnotesize 101}&
{\footnotesize 107}&
&
\tabularnewline
\hline 
{\footnotesize 214}&
{\footnotesize 107}&
{\footnotesize 107}&
&
\tabularnewline
\hline 
{\footnotesize 220}&
{\footnotesize 107}&
{\footnotesize 113}&
&
\tabularnewline
\hline 
{\footnotesize 226}&
{\footnotesize 113}&
{\footnotesize 113}&
&
\tabularnewline
\hline 
{\footnotesize 232}&
{\footnotesize 113}&
{\footnotesize 119}&
{\footnotesize 232}&
{\footnotesize 101+131}\tabularnewline
\hline 
{\footnotesize 238}&
{\footnotesize 119}&
{\footnotesize 119}&
{\footnotesize 238}&
{\footnotesize 107+131}\tabularnewline
\hline 
{\footnotesize 244}&
{\footnotesize 119}&
{\footnotesize 125}&
{\footnotesize 244}&
{\footnotesize 113+131}\tabularnewline
\hline 
{\footnotesize 250}&
{\footnotesize 125}&
{\footnotesize 125}&
{\footnotesize 250}&
{\footnotesize 113+137}\tabularnewline
\hline 
{\footnotesize 256}&
{\footnotesize 125}&
{\footnotesize 131}&
{\footnotesize 256}&
{\footnotesize 107+149}\tabularnewline
\hline 
{\footnotesize 262}&
{\footnotesize 131}&
{\footnotesize 131}&
&
\tabularnewline
\hline 
{\footnotesize 268}&
{\footnotesize 131}&
{\footnotesize 137}&
&
\tabularnewline
\hline 
{\footnotesize 274}&
{\footnotesize 137}&
{\footnotesize 137}&
&
\tabularnewline
\hline
\end{tabular}
\par\end{centering}

\caption{\label{tab:Suma-de-betas,}Possible solution of the problem in two
$\beta$ sums.}
\end{table}

First we can establish that while moving along a number family, $\alpha$
or $\beta$ , the gaps, (the group of compound number inside them)
are placed between prime numbers, which means that, when we find an
even number which is sum of compounds, it would be possible to replace
one element of the sum by the next element of its family, and the
other with the preceding element of its family. This is achieved by
summing 6 in one column and subtracting 6 in the other. Actually,
in Table \ref{tab:Suma-de-betas,} we show even numbers generated
by sums of two $\beta$ . Numbers of each column are glided and each
number is duplicated in order to cover all possibilities. We observe
that it is possible to replace compound sums by prime sums if we replace
the compound factor by the following prime of its family, and the
compound of the other column by the preceding prime of its family.
By doing this way it is possible to find a prime correlation to a
given sum.

\section{About the introduction of a graphic mode for the conjecture}

We will begin to observe the list of numbers $\alpha$ and $\beta$
from a non mathematic point of view, something like a binary series
of white and black boxes. Where black represents compound numbers
and white represents prime numbers. By imaging this kind of bands,
we see that when comparing any band with itself, we are observing
sums of equal factors. On Table \ref{tab:Suma-de-betas,} this happens
one row at a time. To analyze the replacement of addends of all sums
who have at least one compound as addend we will have to compare a
band with itself inverted. With this kind of configuration coincidences
of whites will represent the adequate sum of two primes.

\begin{figure}
\begin{centering}
\includegraphics[width=8cm,keepaspectratio]{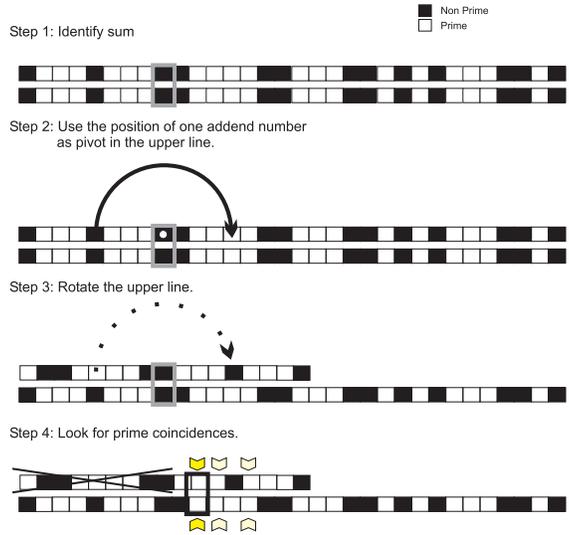}
\par\end{centering}

\caption{\label{Grafo1}Steps to establish the sum relation}
\end{figure}

What we mentioned above can be observed in Figure \ref{Grafo1}, where
black represents compound numbers and white represents prime numbers.
In this graphic we indicate the steps to follow in order to establish
a sum of two equal or different $\alpha$ and the way to express this
sum as a two-prime one.

\begin{figure}
\begin{centering}
\includegraphics[width=8cm,keepaspectratio]{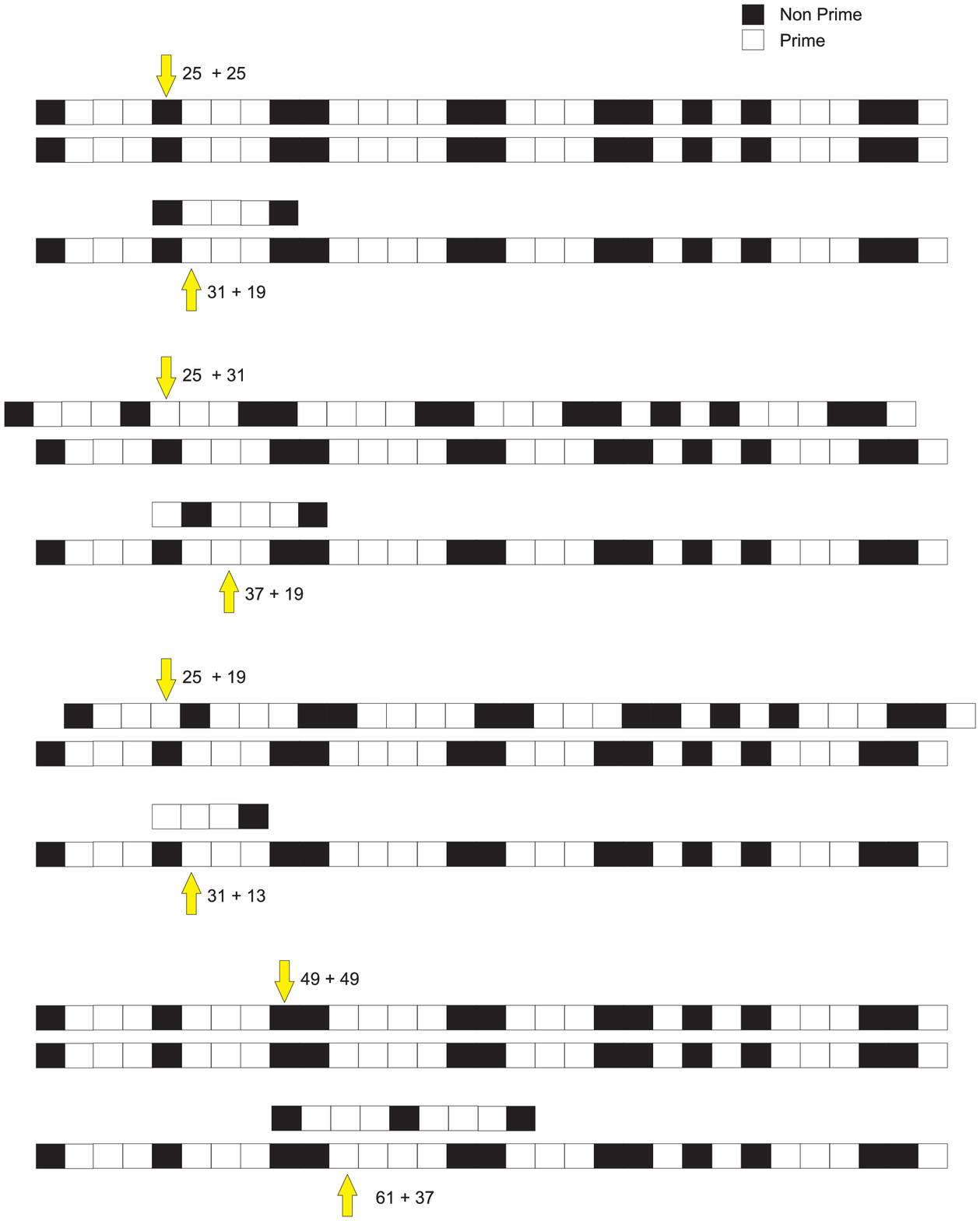}
\par\end{centering}

\caption{\label{Grafo2}Sum relation for the first compound $\alpha$ and
for the second}
\end{figure}

In Figure \ref{Grafo2} we see the sum relation to the following:
$25+25$, $25+31$, $25+19$, and $49+49$

\section{About binary bands}

We will begin our analysis of bands by introducing the idea of band
density. 

In Figure \ref{fig:Densidades-alfa} we can see the prev-density and
the after-density of a band of $\alpha$ numbers, according to the
numbers used as axis to rotate the band. 

\noindent \begin{center}
\begin{figure}
\begin{centering}
\includegraphics[width=10cm,keepaspectratio]{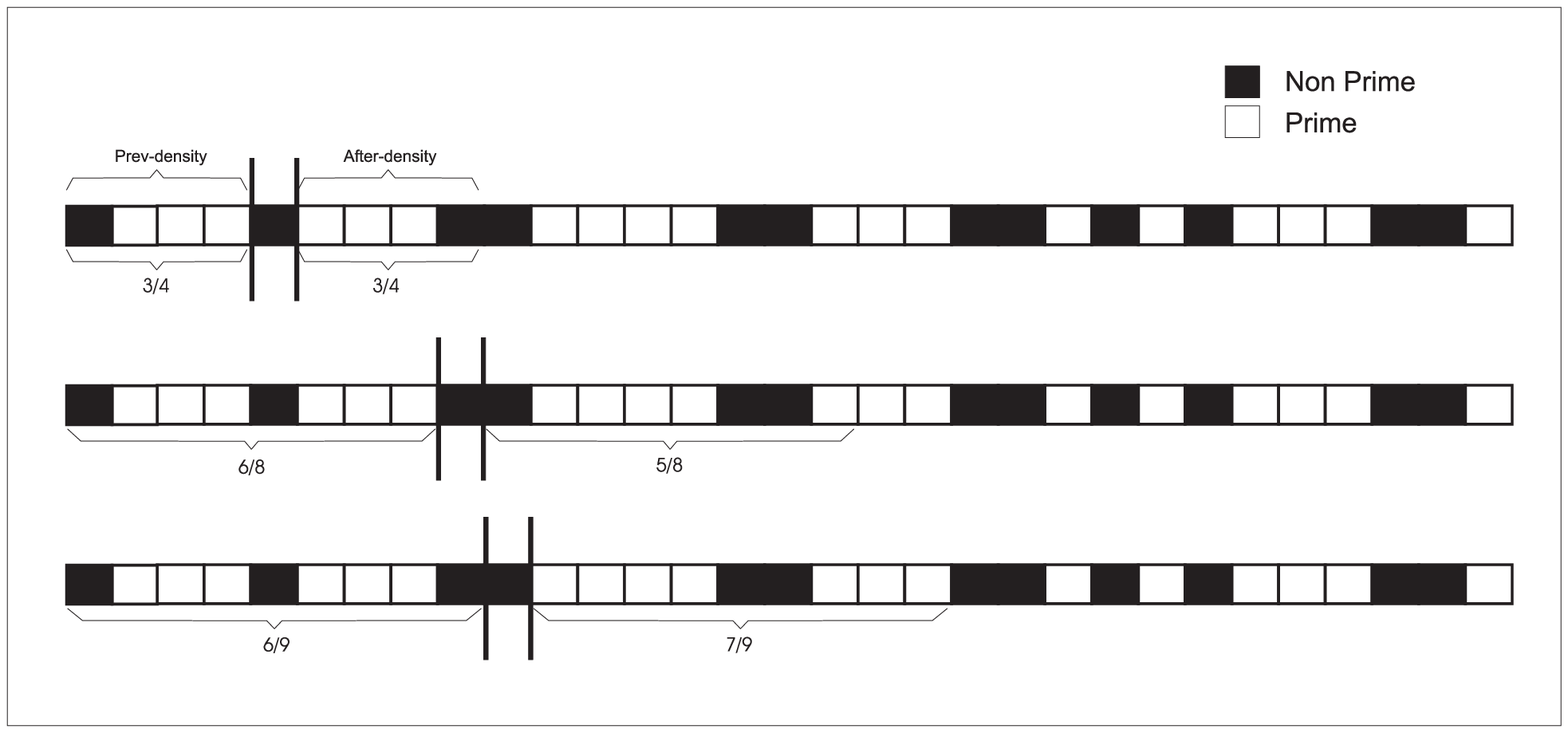}
\par\end{centering}

\caption{\label{fig:Densidades-alfa}Densities for a binary band of $\alpha$
numbers}
\end{figure}

\par\end{center}

Then, we will call \textit{prev-density} to the density of the semi-band
placed previous to the compound number selected as axis. We will call
\textit{after-density} to the density of the semi-band placed after
the axis number, which length is equal to the previous semi-band.

Given a band and inside it a compound number used as axis, we have
three possibilities to fulfill the completeness in order to generate
a family of even numbers:

\begin{enumerate}
\item Glide a band to the right.
\item Glide a band to the left.
\item Compare with itself (No glide needed).
\end{enumerate}
All these possibilities are represented on Figure \ref{fig:Posibles-formas-de}.
After glide is done we have to compare the after-parts of the band
and the prev-parts of it according to the steps indicated on Figure
\ref{Grafo1}.

\begin{center}
\begin{figure}
\begin{centering}
\includegraphics[width=8cm,keepaspectratio]{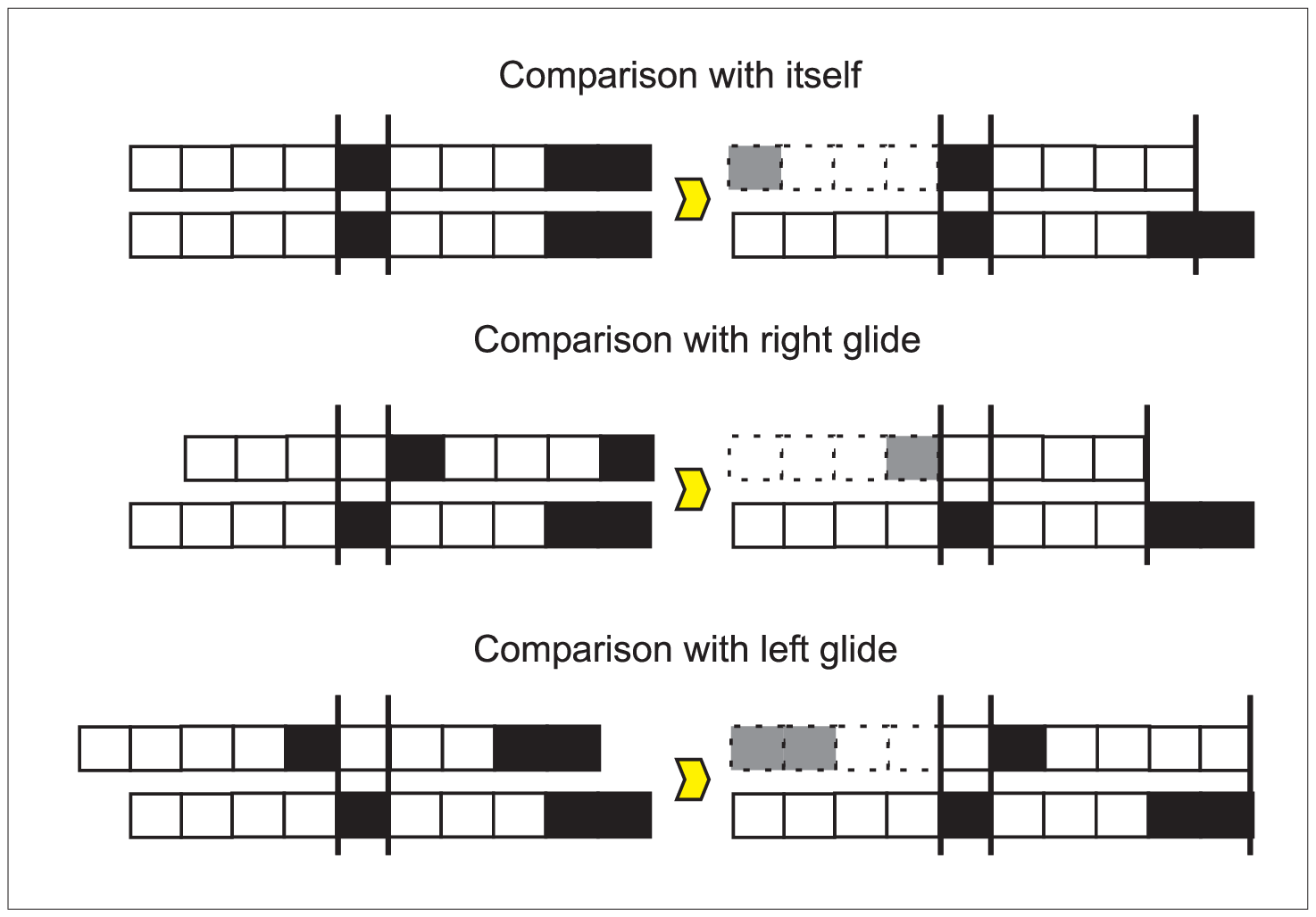}
\par\end{centering}

\caption{\label{fig:Posibles-formas-de}Possible forms of comparing densities}
\end{figure}

\par\end{center}

Using a random band, in which we have left {}``free'' the after-part
for future purposes, our objective will be to obtain a pattern for
which the band could not have any prime coincidence. This pattern
will be called \textit{Non-coincidence pattern}, and it is showed
on Figure \ref{fig:determinaci=F3n-del-patr=F3n}.

\noindent \begin{center}
\begin{figure}
\noindent \begin{centering}
\includegraphics[width=11cm,keepaspectratio]{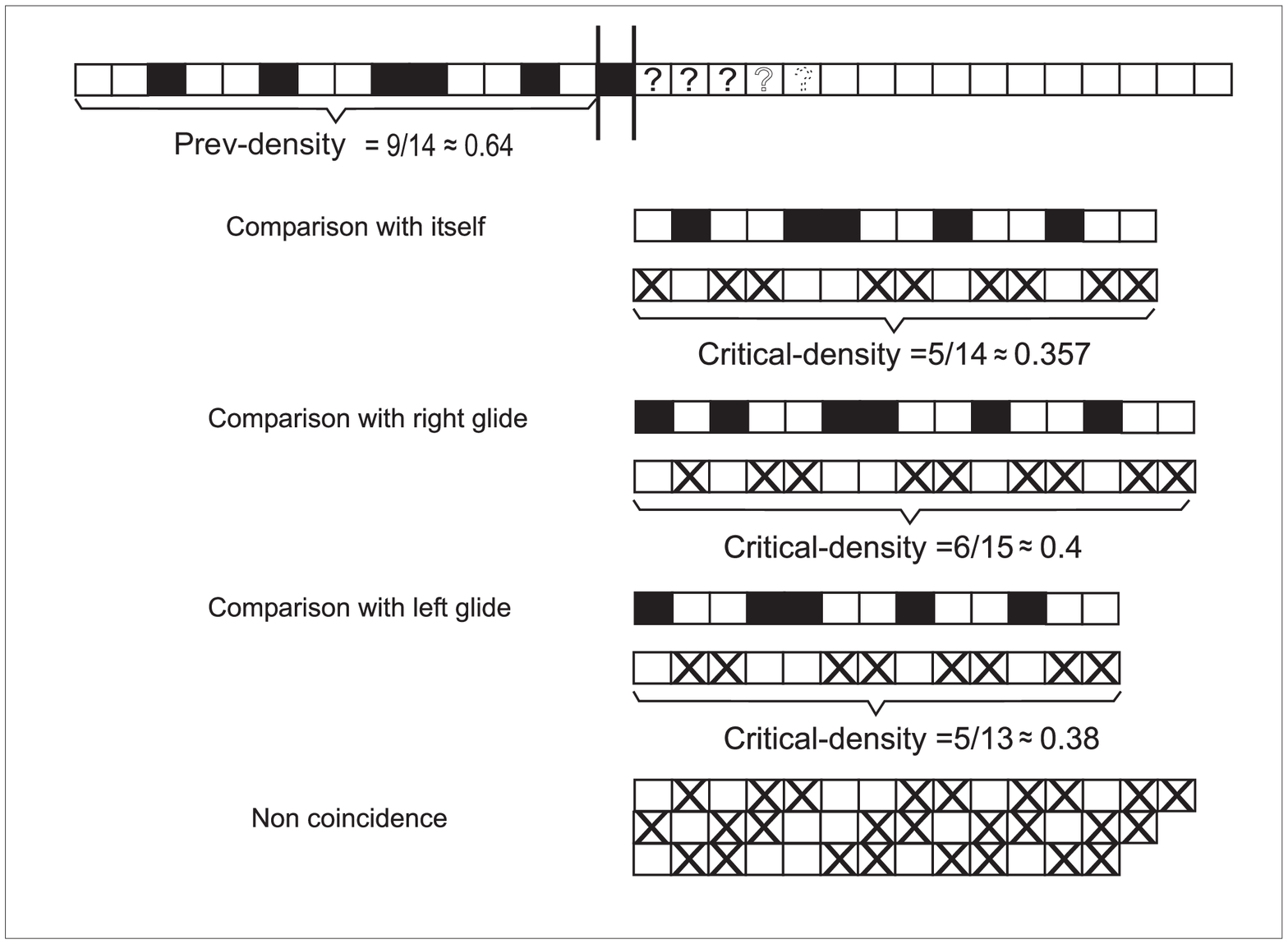}
\par\end{centering}

\caption{\label{fig:determinaci=F3n-del-patr=F3n}Determination of the \textit{non
coincidence pattern}}
\end{figure}

\par\end{center}

We observe that this pattern is the inverse of the prev-band, and
its density correspond to:

\[
1-DA\]

Where $DA$ is the \textit{prev-density}. Since the after-band must
be unique, when superposing all three possible patterns we observe
that the \textit{after-density} should be cero. The density defined
as $1-DA$ will be called \textit{critical-density}.

To observe the behavior of this \textit{critical-density} we will
see the superposition of the after-band of our random band with its
prev-band according to all possible glides. This is shown on Figure
\ref{fig:An=E1lisis-de-la}.

\begin{center}
\begin{figure}
\noindent \begin{centering}
\includegraphics[width=8cm,keepaspectratio]{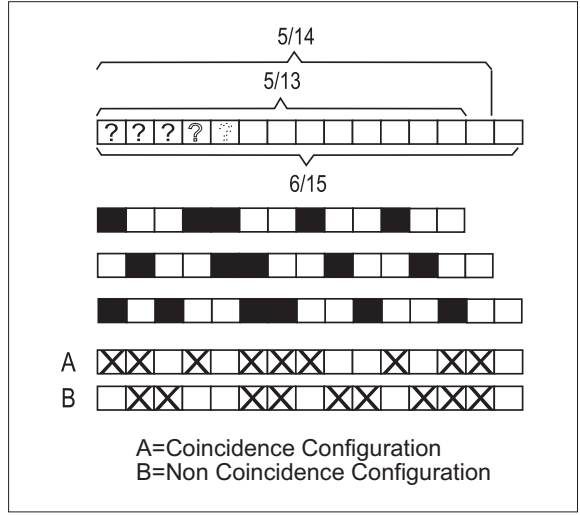}
\par\end{centering}

\caption{\label{fig:An=E1lisis-de-la}Analysis of an after band at \textit{critical-density}}
\end{figure}

\par\end{center}

Suppose we want to place white and black boxes in the after-band according
to the \textit{critical-density}. While placing them in a random manner,
we see that all glides may obtain coincidence. On the other hand,
while keeping a constant density and make coincidence in the last
box of the longest band (left glide), we have only two non coincidence
possibilities for the other glides. It means, placing blacks where
whites are placed for the right glide band, or doing the same for
the other band (no glide one). Since we have kept the \textit{critical-density}
constant, the generated patterns will be the non coincidence ones
at least for one of the bands. When dropping the density down below
the \textit{critical-density} the possibilities of coincidence will
begin to diminish. In this particular case the glided bands are coincident
at least by two, which means that the minimal coincidence density
is of the form $2/\nu$, with a band-length of $\nu$. In other words,
we place a white box in a place of two band coincidence, and the other
white in a coincidence place of the third. In this particular case,
an after-density of $1/\nu$, will lead us to non triple coincidence.
When using bands who have triple coincidence, minimal coincidence
density will be $1/\nu$, and the only chance of non coincidence will
be with after-density equal to cero.

By rising density above \textit{critical}, we see that this growth
must be done by sections, from the longest band to the shortest one.

\noindent \begin{center}
\begin{figure}[H]
\noindent \begin{centering}
\includegraphics[width=8cm,keepaspectratio]{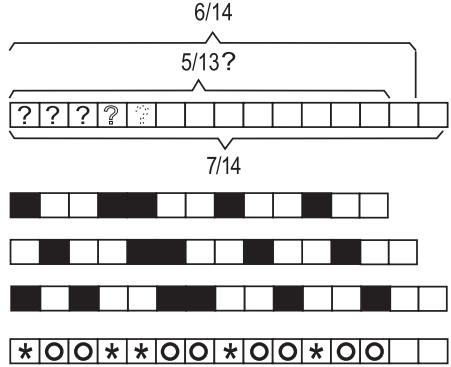}
\par\end{centering}

\caption{\label{fig:Aumento-de-la}Rise of density above critical value}
\end{figure}

\par\end{center}

In Figure \ref{fig:Aumento-de-la} we observe that rising densities
of the longest bands while keeping the density of the shortest (inner
section), we can still place whites and blacks in a non coincidence
pattern with the inner section. This is because we have kept inner
section density equal to critical. Nevertheless, being this the worst
case scenario, we see that by rising density immediately above critical,
the only possibility of placing whites is in a coincidence position.
In Figure \ref{fig:Aumento-de-la} we have placed whites in asterisk
positions according to non coincidence pattern, and boxes with circles
inside represent places where we can place a white box if we rise
density. According to this, we see that when rising density above
critical levels we obtain total coincidence of all glides.

\section{About the use of $\gamma$ family}

We will use $\gamma$ number family (see equation \ref{eq:gamma}),
which contains number $3$ and all its multiples, to make an addition
with the well known families $\alpha$ and $\beta$:

\[
\alpha+\gamma=1+6n_{1}+3+6n_{2}\equiv4+6\eta\]

\[
\beta+\gamma=5+6n_{1}+3+6n_{2}\equiv8+6\chi\equiv2+6\xi\]

\[
\gamma+\gamma=3+6n_{1}+3+6n_{2}\equiv6+6\psi\]

By this way we have generated all families of even numbers, although
$\gamma$ family contains only number $3$ as prime. When using this
number in the conjecture, only sums where $\alpha$ and $\beta$ are
primes are valid, leaving gaps of even numbers that can not be generated
in this fashion. Finally, when adding two $\gamma$, it is clear that
the only even number generated as a sum of two primes is $6$ leaving
behind an entire family of even numbers outside the conjecture. This
justifies why we make use only of the families $\alpha$ and $\beta$.

\section{Conclusions}

The classification of numbers in families -like $\alpha$ and $\beta$-
has been an important factor referring to the organization of numbers.
By the use of such organization, it is possible to follow the path
traced by numbers while adding to each other, and verifying which
family of even numbers have generated. At the same time we have established
that not all even numbers can be written as the sum of any odd prime.
Only numbers of the families $\alpha$ and $\beta$ generate established
even number families.

With binary bands we have presented a new way of watching Goldbach
conjecture from a non mathematic point of view, we hope our analysis
could be the starting point to establish a solution to this problem. 

In future works we expect to perform a detailed study about $\pi$
function. This will allow us a complete analysis of densities of $\alpha$
and $\beta$ numbers, which will have direct application to the band
analysis.

\end{document}